\documentclass[reqno]{amsart}
\usepackage{amsmath,amssymb,amsbsy,amsfonts,amsthm,latexsym,amsopn,amstext,amsxtra,euscript,amscd}
\usepackage{graphicx}
\usepackage{latexsym}
\usepackage{epsfig}

\newtheorem{theorem}{Theorem}
\newtheorem{proposition}{Proposition}
\newtheorem{corollary}{Corollary}

\newtheorem{problem}{Problem}

\def\ds{\displaystyle}
\def\n{\noindent}

\smallskip
\def\bsq{\blacksquare}
\def\eproof{$\hfill\bsq$\par}
\Large
\begin{document}
\title[\bf Moments and the Range of the Derivative]{\bf Moments and the Range of the Derivative }

\author{Eugen J. Ionascu}
\address{Department of Mathematics, Columbus State University, Columbus, GA
31907, US, Honorific Member of the Romanian Institute of
Mathematics ``Simion Stoilow''}
\email{ionascu$\_$eugen@colstate.edu} 

\author{Richard Stephens}
\address{Department of Mathematics, Columbus State University, Columbus, GA
31907, US}

\email{stephens$\_$richard2@colstate.edu}

\date{May $3^{rd}$, 2011}

\begin{abstract} In this note we introduce three problems related
to the topic of finite Hausdorff moments. Generally speaking,
given the first $n+1$ ($n\in \mathbb N\cup \{0\}$) moments,
$\alpha_0$, $\alpha_1$,..., $\alpha_n$, of a real-valued
continuously differentiable function $f$ defined on $[0,1]$, what
can be said about the size of the image of $\frac{df}{dx}$? We
make the questions more precise and we give answers in the cases
of three or fewer moments and in some cases for four moments. In
the general situation of $n+1$ moments, we show that the range of
the derivative should contain the convex hull of a set of $n$
numbers calculated in terms of the Bernstein polynomials,
$x^k(1-x)^{n+1-k}$, $k=1,2,...,n$, which turn out to involve
expressions just in terms of the given moments $\alpha_i$,
$i=0,1,2,...n$. In the end we make some conjectures about what may
be true in terms of the sharpness of the interval range mentioned
before.
\end{abstract} \maketitle

\section{Introduction} We are studying here a problem from real
analysis which can be roughly stated in the following way:

\begin{quote}
\it given a continuously differentiable function whose first $n$
moments are prescribed, what can be said about the image of the
derivative of this function ?
\end{quote}

One of the tools that we will use is the following classical so
called first mean value theorem for integrals (see Section 30.9 in
\cite{bartle}).

\begin{theorem}\label{mean} Let $h$ be a continuous function on $[a,b]$ and
$g$ a non-negative Riemann integrable function. Then there exists
a value $c\in (a,b)$ such that
$$\int_a^bh(x)g(x)dx=h(c)\int_a^bg(x)dx.$$
Moreover,  if $h(x)\ge h(c)$ (or $h(x)\le h(c)$) for all $x\in
[a,b]$, then $h(x)=h(c)$ for every $x$ point of continuity of $g$
and $g(x)>0$.
\end{theorem}

To introduce our hypothesis we let $n\in \mathbb N\cup \{0\}$ and
let $f$ be a continuously differentiable function which satisfies
the following Hausdorff moment type interpolation conditions:

\begin{equation}\label{thegiven}
\int_0^1 x^kf(x)dx=\alpha_k,\ \ \  k=0,1,2,...,n, \ \ \
\alpha_k\in \mathbb R.
\end{equation}

Let us observe that given arbitrary moments $\alpha_k$ the system
(\ref{thegiven}) leads to a linear one if $f$ is a polynomial
function. The main matrix of the resulting system is a Hilbert
matrix. This type of matrix is well know (see \cite{mch}, for
instance) and has a non-zero determinant.

Our investigation was motivated by a proposed problem in the
College Mathematics Journal (\cite{dvt}) which requires one to
show that if $n=2$ and $\alpha_k=k+1$, there exist $c_1, c_2\in
[0,1]$ such that $f'(c_1)=-24$ and $f'(c_2)=60$. It turns out that
this problem was inspired by a problem of C. Lupu (see
\cite{cezarlupu}) which referred to only  two moments,
$\alpha_0=\alpha_1=1$, and asked for a point $c$ where $f'(c)=6$.
We wondered if these numbers were, in a certain sense which will
be defined next, sharp. We will show that this is indeed the case
in the next section (Theorem~\ref{firsttheorem}). Similar
optimization questions, given the first $n$ Hausdorff moments on
$[0,1]$ or $[-1,1]$, are customary subjects in the literature (see
\cite{gi1}, \cite{rodriguez-seatzu}) We are going to formulate the
following very general questions that are our main interest in
this paper.

\begin{problem}\label{problem1} {\it For a fixed $n$ and $\alpha_k$ as before, what is the largest range $[A,B]$
such that  $[A,B]\subseteq Range(f')$ for every $f$ a continuously
differentiable function on $[0,1]$ satisfying (\ref{thegiven})?}
\end{problem}

\begin{problem}\label{problem2} {\it For a fixed $n$ and $\alpha_k$ as before, what is the biggest number $L$ such
that for every $f$ a continuously differentiable function on
$[0,1]$ satisfying (\ref{thegiven}) there exists some interval
$[a,b]$ with $b-a=L$ that satisfies $[a,b]\subseteq Range(f')$ ?}
\end{problem}

We observe that in order to prove that $[A,B]$ is the answer for
Problem~\ref{problem1},  it is necessary to show that
$[A,B]\subseteq Range(f')$ for every $f$ a continuously
differentiable function on $[0,1]$ satisfying (\ref{thegiven}) and
that for every $\epsilon>0$ there exists $f_l$ and $f_r$
continuously differentiable functions on $[0,1]$ satisfying
(\ref{thegiven}) and

\begin{equation}\label{equation1} Range(f_l')\subseteq (A-\epsilon,\infty)\
 and \ Range(f_r')\subseteq (-\infty,B+\epsilon).
 \end{equation}

It is clear that if $A$ and $B$ give the answer in
Problem~\ref{problem1}, then in trying to answer
Problem~\ref{problem2} we must have $L\ge B-A$. If for every
$\epsilon>0$, one can find a function ($f_l=f_r$) that satisfies
both conditions in (\ref{equation1}), then the answer to
Problem~\ref{problem2} is simply $L=B-A$.

Another related problem here is to characterize the case $L> B-A$
and calculate $L$ in this case in terms of the $\alpha_k$'s.
Perhaps Problem 2 may be easier if one restricts the class of
functions in consideration to something more manageable like
polynomials of a certain degree.

If we want to make the range of the derivative as small as
possible, we just have to take moments that satisfy the necessary
and sufficient condition  for having a solution to the  system
that results from having a linear function, say $f(x)=u+vx$, $x\in
[0,1]$, satisfying (\ref{thegiven}):

$$\frac{u}{k+1}+\frac{v}{k+2}=\alpha_k,\ \ k=0,1,2,....n.$$

This is equivalent to

$$rank \left[%
\begin{array}{ccc}
  1 & \frac{1}{2} & \alpha_0 \\
  \frac{1}{2} & \frac{1}{3} & \alpha_1 \\
  \frac{1}{3} & \frac{1}{4} & \alpha_2 \\
  ... & & \\
  \frac{1}{n} & \frac{1}{n+1} & \alpha_{n-1}\\
  \frac{1}{n+1} & \frac{1}{n+2} & \alpha_{n}
\end{array}%
\right]=2.$$

On the other hand, if we want to make the range of $f'$ as big as
possible, it makes sense to restrict our moments to a finite
range, say $[-1,1]$. We observe that the problem is homogeneous
under dilations, so let us formulate a third problem.

\begin{problem}\label{problem3} {\it For a fixed $n$,
what is the maximum of  $B-A$ such that  $[A,B]\subseteq
Range(f')$ for every $f$ a continuously differentiable function on
$[0,1]$ satisfying (\ref{thegiven}), the maximum being taken over
all possible moments $\alpha_k\in [-1,1]$?}
\end{problem}

\n We will show in Section 2 that the answer to
Problem~\ref{problem3} is 156 if $n=2$, for the moments
$\alpha_0=1$, $\alpha_1=-1$ and $\alpha_2=1$. We observe that if
the answer to Problem~\ref{problem2} is zero, then the answer to
Problem~\ref{problem3} is also zero.  As suggested by one of the
referees of our paper, one can ask similar questions about the
range of $f''$ or higher derivatives, assuming these exist. We
will make some remarks about these questions and see how the
results for the first derivative could be applied for higher
derivatives.

\vspace{0.2in}
\section{Small values of $n$}
We have a few complete answers to Problem~\ref{problem1} for small
values of $n$ ($n\le 3$).  First, let us study what happens with
$n=0$. If we take $g(x)=1-x$ and $h=f'$ in Theorem~\ref{mean},
using integration by parts, we get
$$f'(c_1)\frac{1}{2}=\int_0^1f'(x)(1-x)dx=-f(0)-\int_0^1f(x)(-1)dx=\alpha_0-f(0), \ c_1\in (0,1),$$
\n or
$$f'(c_1)=2(\alpha_0-f(0)), \ c_1\in (0,1).$$

\n If $f(0)=a$, then we can take $f(x)=a+(2\alpha_0-2a)x$ and
observe that in case $n=0$, there exists a function such that
$\int_0^1f(x)dx=\alpha_0$ and $Range(f')=\{2\alpha_0-2f(0)\}$.
This gives us the following simple answers to
Problem~\ref{problem1} and Problem~\ref{problem2}.

\begin{proposition}\label{firstprop} For $n=0$, there is no $A$
and $B$ that satisfy the requirements of Problem~\ref{problem1}.
The answer for Problem~\ref{problem2} ($n=0$) is $L=0$.
\end{proposition}

Let us continue the analysis in the case $n=1$. We can apply
Theorem~\ref{mean} to $g(x)=x(1-x)$ and $h=f'$, $x\in [0,1]$.
Then, a similar  calculation gives that for some $c_2\in (0,1)$,

$$f'(c_2)\frac{1}{6}=-\int_0^1f(x)(1-2x)dx=2\alpha_1-\alpha_0,   \Rightarrow f'(c_2)=6(2\alpha_1-\alpha_0).$$

If we apply Theorem~\ref{mean} to $h=f'$ and $g(x)=(1-x)^2$
instead,

$$f'(c_3)\frac{1}{3}=\int_0^1f'(x)(1-x)^2dx=-f(0)-\int_0^1f(x)(2x-2)dx=2(\alpha_0-\alpha_1)-f(0), \ c_3\in (0,1),$$
\n or
$$f'(c_3)=6(\alpha_0-\alpha_1)-3f(0),\ \text{for \ some}\ \ c_3\in (0,1).$$ So,
if we take $a=2(2\alpha_0-3\alpha_1)$ and $f(x)=a+mx$ where
$m=6(2\alpha_1-\alpha_0)=6(\alpha_0-\alpha_1)-3a$, we get a
function which will give us what we need in this case, and
therefore provide a similar answers to our problems.

\begin{proposition}\label{secondprop} For $n=1$,  we can take $A=B=12\alpha_1-6\alpha_0$
to satisfy the requirements of Problem~\ref{problem1}. The answer
for Problem~\ref{problem2} ($n=1$) is $L=0$.
\end{proposition}

The case $n=2$ is getting a little more interesting; it is
essentially non-trivial and at the same time pretty surprising. We
have a definite  answer to Problem~\ref{problem1} and
Problem~\ref{problem3} and we show some inequality for $L$ in
Problem~\ref{problem2}.

\begin{theorem}\label{firsttheorem} For $n=2$,  if \ $\Delta_0:=6\alpha_2-6\alpha_1+\alpha_0>0$,
the values

$$A:=12(4\alpha_1-\alpha_0-3\alpha_2)\ \text{and}\ B:=12(3\alpha_2-2\alpha_1)$$
 satisfy the requirements of
Problem~\ref{problem1} and if $\Delta_0<0$ then one needs to
switch the values of $A$ and $B$ above in order to solve
Problem~\ref{problem1}. If $\Delta_0=0$, the values
$A=B=12(3\alpha_2-2\alpha_1)$ answer Problem~\ref{problem1} and
$L=0$ answers Problem~\ref{problem2}.
\end{theorem}

\proof First, let us show that $A$ and $B$ are always in the range
of the derivative. This is done as we have seen before by setting
in the Theorem~\ref{mean}, $h=f'$, and $g(x)=x(1-x)^2\ge 0$ ($x\in
[0,1]$). Indeed, we have $\int_0^1g(x)dx=\frac{1}{12}$ and

$$\int_0^1f'(x)g(x)dx=f(x)g(x)|_0^1-\int_0^1f(x)(1-4x+3x^2)dx=4\alpha_1-\alpha_0-3\alpha_2.$$

Hence, for some $c_4$ we must have
$f'(c_4)\int_0^1g(x)dx=4\alpha_1-\alpha_0-3\alpha_2$ which in turn
gives $f'(c_4)=12(4\alpha_1-\alpha_0-3\alpha_2)=A$. Similarly, for
$g(x)=x^2(1-x)$, ($x\in [0,1]$), one finds that
$\int_0^1g(x)dx=\frac{1}{12}$ still holds true and
$$\int_0^1f'(x)g(x)dx=\int_0^1f(x)(2x-3x^2)dx=3\alpha_2-2\alpha_1.$$

This insures that $B=12(3\alpha_2-2\alpha_1)$ is also in the range
of $f'$. Because $f'$ is assumed to be continuous we get that the
whole interval $[A,B]$ or $[B,A]$ is contained in the range of
$f'$.

From here on, we are going to work under the first assumption
($\Delta_0>0$) which is equivalent to $A<B$ ($B-A=12\Delta_0$). To
show that $A$ and $B$ are sharp bounds we begin with $B$ by
constructing a spline function $s_t$ for $t\in (0,1)$, defined by

$$s_t(x)=\begin{cases} a+bx+cx^2\ for \ x\in [0,t]\\ \\
m+nx\ \ if\ \ x\in [t,1],
\end{cases}$$

\n where $a$, $b$,  $c$, $m$ and $n$ are determined by the
conditions $\int_0^1s_t(x)dx=\alpha_0$,
$\int_0^1xs_t(x)dx=\alpha_1 $, $\int_0^1x^2s_t(x)dx=\alpha_2 $ and
the restrictions necessary to insure that $s_t$ is continuously
differentiable at $x=t$.

In order to add an intuition element we included here the graphs
of $s_{1/10}$ and its derivative for $\alpha_0=1$, $\alpha_1=1$
and $\alpha_2=2$.

\begin{figure}
\[
\underset{(a)\ s_{\frac{1}{10}}}
{\epsfig{file=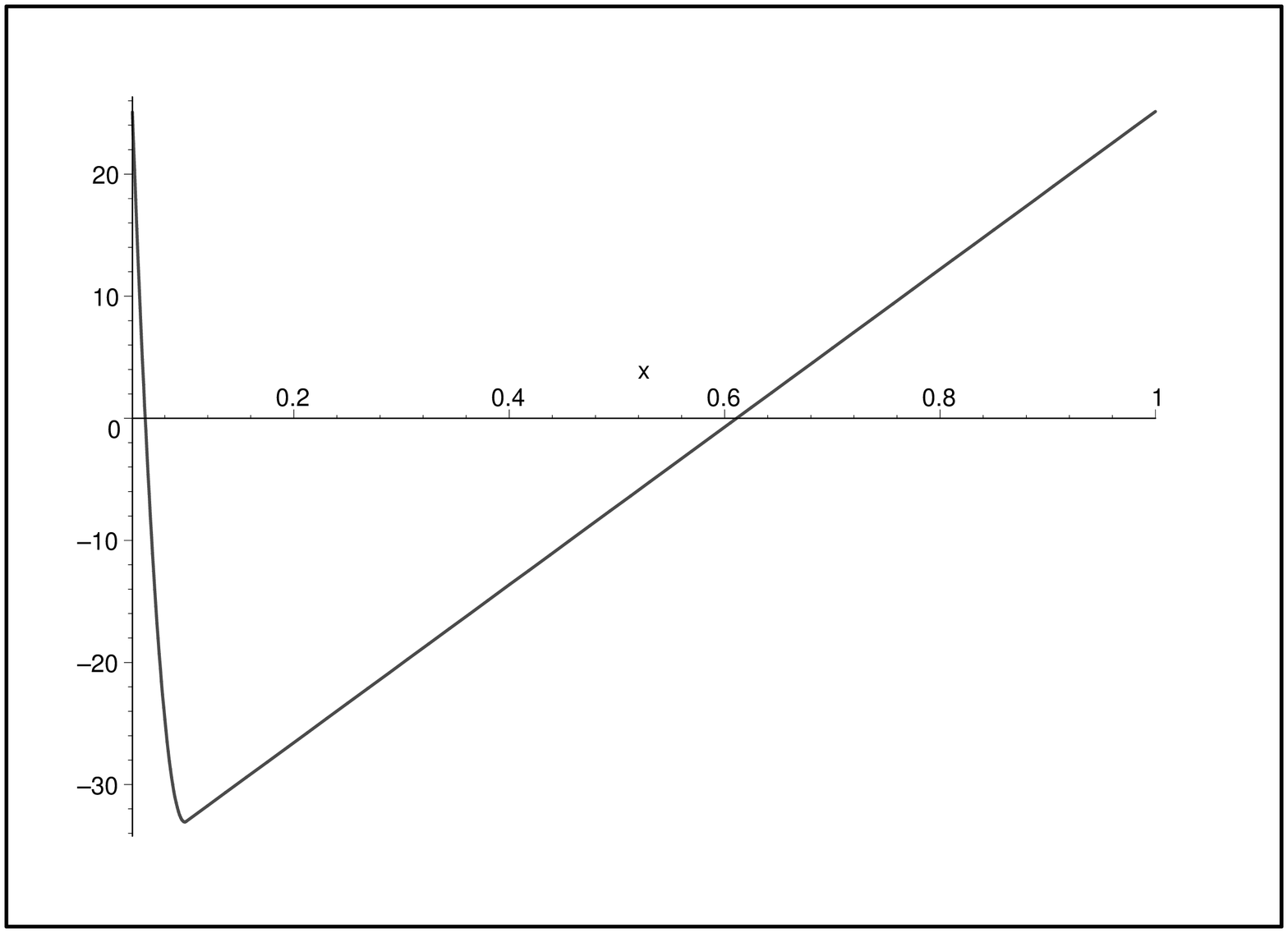,height=2in,width=2in,angle=-90}}\ \ \ \ \ \
\underset{(b) \ s'_{\frac{1}{10}}}
{\epsfig{file=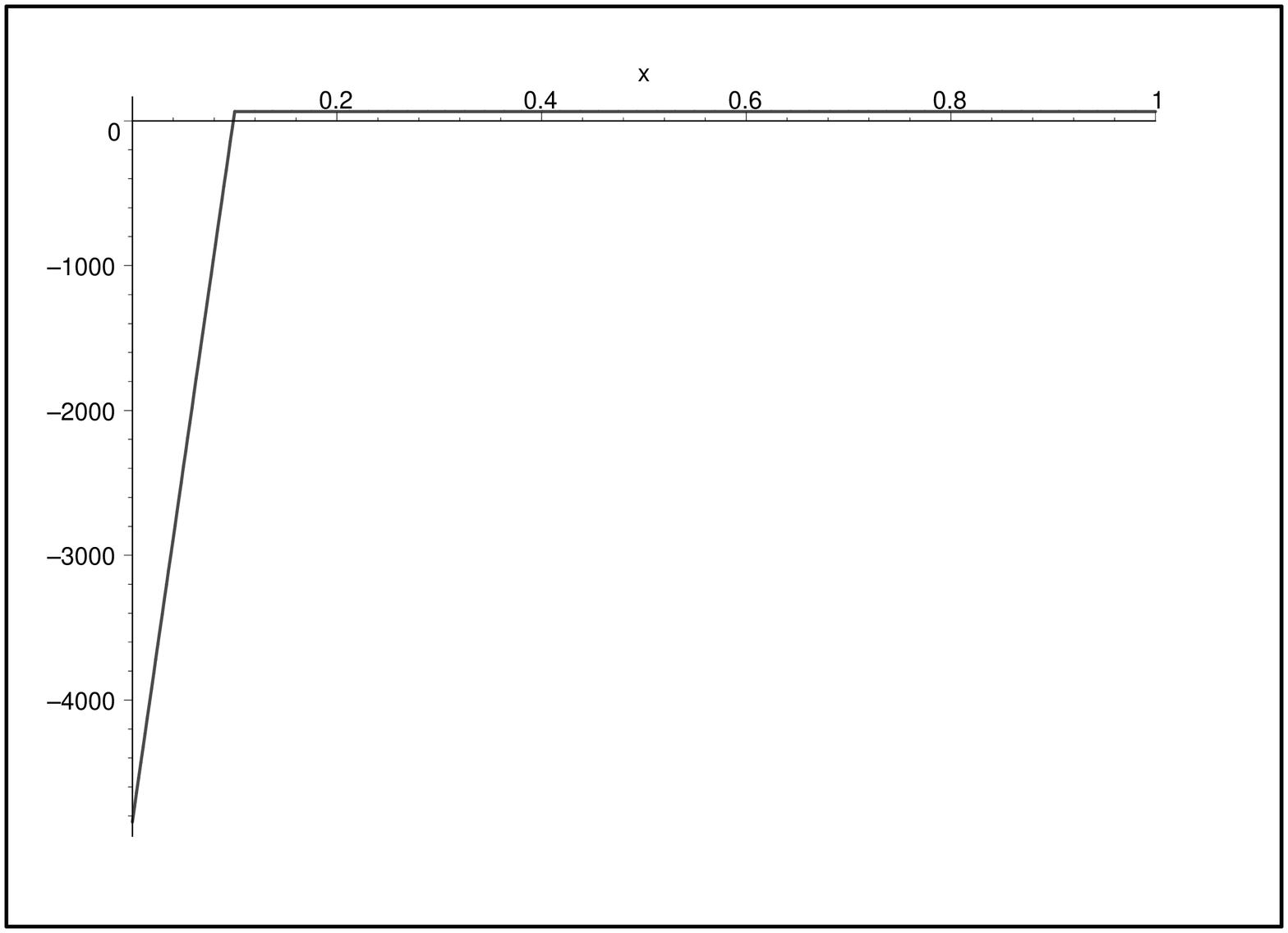,height=2in,width=2in,angle=-90}}\
\]
\caption{Case $n=2$ and $\alpha_k=k+1$}
\end{figure}

It is easy to see that $s_t$ is continuously differentiable at
$x=t$ if and only if $m=a-ct^2$ and $n=b+2ct$. This gives the new
expression of $s_t$ just in terms of $a$, $b$ and $c$:

\begin{equation}\label{equationofspline}
s_t(x)=\begin{cases} a+bx+cx^2\ for \ x\in [0,t]\\ \\
(a-ct^2)+(b+2ct)x\ \ if\ \ x\in [t,1].
\end{cases}
\end{equation}

\n One can check that the moment restrictions reduce to the
following $3\times 3$ relatively simple linear system of equations
in $a$, $b$ and $c$:

\begin{equation}\label{systemforabc}
\begin{cases} \ds a+\frac{b}{2}+\left(\frac{t^3}{3}-t^2+t\right)c=\alpha_0\\ \\
\ds
\frac{a}{2}+\frac{b}{3}+\left(\frac{t^4}{12}+\frac{2t}{3}-\frac{t^2}{2}\right)c
=\alpha_1\\ \\
\ds
\frac{a}{3}+\frac{b}{4}+\left(\frac{t^5}{30}+\frac{t}{2}-\frac{t^2}{3}\right)c
=\alpha_2.
\end{cases}
\end{equation}

\n It is clear that we have a unique solution for this system at
least for infinitely many values of $t$ since the main determinant
of the system is a polynomial in $t$ of degree at most $5$. Let us
observe that

\begin{equation}\label{derivativeofspline}
\frac{d}{dx}(s_t)(x)=\begin{cases} b+2cx\ for \ x\in [0,t]\\ \\
b+2ct\ \ if\ \ x\in [t,1].
\end{cases}
\end{equation}

\n We observe that if $c>0$, the maximum of this function is
$b+2ct$. With a little work one solves the system
(\ref{systemforabc}) and finds that

$$c=\frac{30\Delta_0}{t^3(6t^2-15t+10)},
b+2ct=\frac{12(6t^2\alpha_1-3t^2\alpha_0-20\alpha_1+30\alpha_2+5t\alpha_0-15t\alpha_2)}{6t^2-15t+10}.$$

\n It is clear from these expressions that under our hypothesis
$c>0$ for every $t>0$ and that $\ds
B=12(3\alpha_2-2\alpha_1)=\lim_{t\to 0} (b+2ct)$ which proves that
$B$ is sharp. In a similar way one can prove that $A$ is sharp by
taking a spline $\tilde{s}_t$ which is first a linear piece on
$[0,t]$ and a quadratic piece on $[t,1]$. It turns out that the
calculations are very much similar to the ones above with the only
difference that this time we let $t$ approach $1$.

However, we are going to show that the lower bound $A$ is sharp by
using an invariance principle here by doing a ``change of
variable" so to speak and considering how the
Problem~\ref{problem1} changes from $f$ to $g$ where
$g(x)=f(1-x)$, $x\in [0,1]$. Let us denote by $\alpha_k(f)$ the
$k^{th}$ moment for the function $f$. One can see that

$$\alpha_0(g)=\alpha_0(f),\ \alpha_1(g)=\alpha_0(f)-\alpha_1(f)\ and\ \alpha_2(g)=\alpha_0(f)-2\alpha_1(f)+\alpha_2(f),$$

\n and of course, the relations are symmetric with respect to
interchanging $f$ and $g$ , i.e.

$$\alpha_0(f)=\alpha_0(g),\ \alpha_1(f)=\alpha_0(g)-\alpha_1(g)\ and\ \alpha_2(f)=\alpha_0(g)-2\alpha_1(g)+\alpha_2(g).$$

\n Let us observe that the hypothesis that $\Delta_0(f)>0$ is in
fact invariant under this change:
$$6\alpha_2(f)-6\alpha_1(f)+\alpha_0(f)=6\alpha_2(g)-6\alpha_1(g)+\alpha_0(g)>0.$$

By the first part of our proof, we see that
$$B(g)=12(3\alpha_2(g)-2\alpha_1(g))=12[3(\alpha_0(f)-2\alpha_1(f)+\alpha_2(f))-2(\alpha_0(f)-\alpha_1(f)),$$

\n or $$B(g)=12(3\alpha_2(f)-4\alpha_1(f)+\alpha_0(f),$$

\n is a sharp bound for the range of $g'$.  Since $g'(x)=-f'(1-x)$
we see that the range of $g$ is just the range of $f$ reflected
into the origin and vice versa. Hence,
$A(f)=-B(g)=12[4\alpha_1(f)-3\alpha_2(f)-\alpha_0(f)]$ is a sharp
lower bound for $f$. The rest of the statements of the theorem
follow from what we have shown so far.\eproof

\vspace{0.2in}

\begin{corollary}\label{aboutL}
In the case $n=2$, in respect to Problem~\ref{problem2}, we have
$$12|\Delta_0|\le L\le 32|\Delta_0|.$$
The maximum required in Problem~\ref{problem3} is 156.
\end{corollary}

\proof The first part is a simple consequence of the fact
$B-A=12\Delta_0$  and the last part follows from the fact that
$|\Delta_0|=|6\alpha_2-6\alpha_1+\alpha_0|\le 13$ if $\alpha_i\in
[-1,1]$, $i=0,1,2$. To show the inequality $L\le 32|\Delta_0|$ we
employ the same idea by constructing spline which is symmetric
around $1/2$:

\begin{equation}\label{equationofsplinesymmetric}
\hat{s}_t(x)=\begin{cases}  a+bx+c(x-2tx-t^2+t-1/4)\ for \ x\in [0,\frac{1}{2}-t]\\ \\
a+bx+cx^2\ \ if\ \ x\in [\frac{1}{2}-t,\frac{1}{2}+t]\\ \\
a+bx+c(x+2tx-t^2-t-1/4).
\end{cases}
\end{equation}
This spline is continuously differentiable on $[0,1]$ and depends
on three parameters which if determined from the constraints given
by the moments we get $$c=\frac{120 \Delta_0}{t(15 - 40t  +
48t^2)}>0, \text{for all} \ t\in [0,1], $$ and which shows that
the minimum and the maximum of the derivative of $\hat{s}_t$ is
attained on the linear pieces. One can see that the difference
between these two values is actually $4ct$ and so letting $t\to 0$
we get that $$L\le \lim_{t\to 0}4ct=32|\Delta_0|.$$ \eproof

One can use the same techniques to show that for three moments,
assuming the second derivative exists, the range of the second
derivative should contain $30\delta_0$ and this is sharp because a
polynomial of degree two exists solving the moments problem.

The case $n=3$ is even more interesting and a lot more
complicated. First of all we have at least three new possible
values that we need to add to the range of $f'$:

\begin{equation}\label{defofCDE}
 C=20(4\alpha_3-3\alpha_2),\
D=60(3\alpha_2-2\alpha_3-\alpha_1),\ \text{and}\
E=20(4\alpha_3-9\alpha_2+6\alpha_1-\alpha_0)
\end{equation}

\n obtained from Bernstein polynomials, $g_1(x)=x^3(1-x)$,
$g_2(x)=x^2(1-x)^2$ and $g_3(x)=x(1-x)^3$ respectively.

\begin{proposition}\label{enlargement} Given $A$ and $B$ as defined in Theorem~\ref{firsttheorem},
we have the inclusion
$$[min(A,B),max(A,B)]\subset [min(D,C,E),max(D,C,E)].$$
\end{proposition}

\proof Let us observe that $g_1(x)+g_2(x)=x^2(1-x)$ and
$g_2(x)+g_3(x)=x(1-x)^2$. Differentiating and integrating against
$f(x)$ we get the relations
$\frac{1}{30}D+\frac{1}{20}E=\frac{1}{12}A$ or
$A=\frac{2}{5}D+\frac{3}{5}E$ and similarly
$B=\frac{3}{5}C+\frac{2}{5}D$. These two convex linear
combinations are enough to conclude the desired statement. \eproof

Of course, this proposition can be generalized to an arbitrary
$n$. So we expect that the interval that answers
Problem~\ref{problem1} contains the convex hull of the numbers
constructed as usual, i.e.

\begin{equation}\label{definitionOfD}
D_k:=-\frac{\int_0^1f(x)\frac{d}{dx}(x^k(1-x)^{n+1-k})dx}{\int_0^1x^k(1-x)^{n+1-k}dx},
k=1,2,...,n,
\end{equation}

\n given by the highest degree Bernstein basis polynomials
possible.

We observe that if we define
$\Delta_1:=10\alpha_3-12\alpha_2+3\alpha_1$, then
$$C=E+20\Delta_0\ \text{and}
\ D=C-20\Delta_1.$$

\n Hence, we observe that if we have $\Delta_0>0$ and $\Delta_1<0$
for instance, then $D>C>E$. Therefore, in light of
Proposition~\ref{enlargement}, the candidates for the two values
needed to answer Problem~\ref{problem1} are $\tilde{A}=E$ and
$\tilde{B}=D$ under the given assumption. In fact, for various
other situations we believe that the values $A$ and $B$ that
answer Problem~\ref{problem1} are given for each case in the
following table

\vspace{0.1in}

\begin{tabular}{|c|c|c|c|}
  \hline
 No& Hypothesis  & A & B \\
  \hline
 (i) &$\Delta_0\ge 0, \Delta_1\le 0$ & $20(4\alpha_3-9\alpha_2+6\alpha_1-\alpha_0)$ & $60(3\alpha_2-2\alpha_3-\alpha_1)$ \\
 (ii)& $ 0\le \Delta_0 \le \Delta_1$ & $60(3\alpha_2-2\alpha_3-\alpha_1)$ & $20(4\alpha_3-3\alpha_2)$ \\
  (iii) &$ 0\le \Delta_1 \le \Delta_0$ & $20(4\alpha_3-9\alpha_2+6\alpha_1-\alpha_0)$ & $20(4\alpha_3-3\alpha_2)$ \\
  (iv)  &$\Delta_0\le  0, \Delta_1\ge 0$ &  $60(3\alpha_2-2\alpha_3-\alpha_1)$ &$20(4\alpha_3-9\alpha_2+6\alpha_1-\alpha_0)$\\
 (v)& $ \Delta_1 \le \Delta_0\le 0$ &  $20(4\alpha_3-3\alpha_2)$ &$60(3\alpha_2-2\alpha_3-\alpha_1)$ \\
  (vi) &$ \Delta_0 \le \Delta_1\le 0$ &  $20(4\alpha_3-3\alpha_2)$ & $20(4\alpha_3-9\alpha_2+6\alpha_1-\alpha_0)$ \\
  \hline
\end{tabular}\par

\vspace{0.1in}
\n where $\Delta_0=72 \left|%
\begin{array}{ccc}
   1 & \frac{1}{2} & \alpha_0 \\
  \frac{1}{2} & \frac{1}{3} & \alpha_1 \\
  \frac{1}{3} & \frac{1}{4} & \alpha_2 \\
\end{array}%
\right|$ and
$\Delta_1=720\left|%
\begin{array}{ccc}
  \frac{1}{2} & \frac{1}{3} & \alpha_1 \\
  \frac{1}{3} & \frac{1}{4} & \alpha_2 \\
   \frac{1}{4} & \frac{1}{5} & \alpha_3 \\
\end{array}%
\right|.$

We have the following partial result along these lines.

\begin{theorem}\label{secondtheorem} For $n=3$, the upper bound of
(i) and the lower bound of (ii), in the table above, are correct.
 If $\Delta_0=\Delta_1=0$ then $A=B=20(4\alpha_3-3\alpha_2)$ and $L=0$
solves Problem~\ref{problem1} and Problem~\ref{problem2}.
\end{theorem}

\proof First of all let us observe that the cases (iv), (v) and
(vi) follow from (i), (ii) and  (iii) respectively by simply
changing $f$ into $-f$. This simple transformation changes
basically the order of $A$ and $B$. It is easy to see that
$\Delta_0=\Delta_1=0$ implies the existence of a linear map that
has the given moments and so $A=B$ and $L=0$. Hence in what
follows we will assume that $\Delta_0\not = 0$ or $\Delta_1\not
=0$.

Based on the invariance principle that we used in the proof of
Theorem~\ref{firsttheorem} we need to show the sharpness of only
the upper bound in (i). Indeed we observe that if $g(x)=f(1-x)$,
$x\in [0,1]$ then one can check that the hypothesis
$\Delta_1(f)\le 0$ changes into $\Delta_0(g)\le \Delta_1(g)$.
Also, the hypothesis $0\le \Delta_1 \le \Delta_0$, is actually
invariant under this change. One also needs to take into account
that the bound $D$ is invariant under this transformation but $C$
and $E$ interchange:

$$D(f)=D(g),\ C(f)=E(g), \ \text{and}\ E(f)=C(g).$$

So, let us begin with case (i) and show that
$B=60(3\alpha_2-2\alpha_3-\alpha_1)$ is sharp. For every $t\in
(0,1/2)$, consider a spline function $s_{1,t}$ which is quadratic
on $[0,t]$, linear on $[t,1-t]$ and another quadratic on
$[1-t,1]$. The constraints of having this spline a continuous and
differentiable function give us a similar form for $s_{1,t}$ to
the one constructed in the proof of Theorem~\ref{firsttheorem} in
equality (\ref{equationofspline}), in terms of four free
parameters $a$, $b$, $c$ and $d$:

$$
s_{1,t}(x)=\begin{cases} a+bx+cx^2\ for \ x\in [0,t]\\ \\
(a-ct^2)+(b+2ct)x\ \ if\ \ x\in [t,1-t]\\ \\
d(1-t)^2+a-ct^2+[b+2ct-2d(1-t)]x+dx^2 \ \ if\ \ x\in [1-t,1].
\end{cases}
$$

\n The four parameters are then determined by imposing the four
linear constraints given by the moments. The resulting system is

\begin{equation}\label{systemforabcd}
\begin{cases} \ds a+\frac{b}{2}+(\frac{t^3}{3}-t^2+t)c+\frac{t^3}{3}d =\alpha_0\\ \\
\ds
\frac{a}{2}+\frac{b}{3}+\left(\frac{t^4}{12}+\frac{2t}{3}-\frac{t^2}{2}\right)c+\left(\frac{t^3}{3}-\frac{t^4}{4}\right)d
=\alpha_1\\ \\
\ds
\frac{a}{3}+\frac{b}{4}+\left(\frac{t}{2}-\frac{t^2}{3}+\frac{t^5}{30}\right)c+\left(\frac{t^3}{3}-\frac{t^4}{6}+\frac{t^5}{30}\right)d
=\alpha_2\\ \\
\ds
\frac{a}{4}+\frac{b}{5}+\left(\frac{2t}{5}-\frac{t^2}{4}+\frac{t^6}{60}\right)c+\left(\frac{t^3}{3}-\frac{t^4}{4}+\frac{t^5}{10}-\frac{t^6}{60}\right)d
=\alpha_3.
\end{cases}
\end{equation}

As we have observed before, the system has a unique solution for
infinitely many values of $t\in (0,1/2)$, since the main
determinant of the system is a polynomial in $t$ of degree at most
11.  Because the derivative of $s_{1,t}$ is given by

$$
s'_{1,t}(x)=\begin{cases} b(t)+2c(t)x\ for \ x\in [0,t]\\ \\
b(t)+2c(t)t\ \ if\ \ x\in [t,1-t]\\ \\
b(t)+2c(t)t-2d(t)(1-t-x) \ \ if\ \ x\in [1-t,1].
\end{cases}
$$

\n One can use a symbolic calculator and check that $$\lim_{t\to
0} b(t)+2c(t)t=60(3\alpha_2-2\alpha_3-\alpha_1),$$ \n which is one
necessary fact to prove the sharpness of $B$. Also, we need to
check that for most of the values of $t$, $b(t)+2c(t)t$ is a
maximum of the derivative of $s_{1,t}$. For this end, it is enough
to check that $c(t)>0$ and $d(t)<0$ for small values of $t$.
Again, one can compute $\lim_{t\to 0}c(t)t^3=\Delta_0-\Delta_1> 0$
under our assumption in case (i) (unless both numbers $\Delta_0$,
$\Delta_1$ are zero). Also, the limit of $d(t)t^3$ as $t\to 0$
turns out to be equal to $3\Delta_1\le 0$. If $\Delta_1 =0$ we
know that $\Delta_0>0$. In this case we have $\lim_{t\to
0}d(t)t^2=-\frac{9}{4}\Delta_0 <0$.

Using the duality via $g(x)=f(1-x)$, we see that
$A=60(3\alpha_2-2\alpha_3-\alpha_1)$ is a sharp lower bound in the
case (ii). \eproof

In the case $n=3$, assuming the table before
Theorem~\ref{secondtheorem} is correct, with respect to
Problem~\ref{problem2}, we have either $L\ge 20|\Delta_0|$, $L\ge
20|\Delta_1|$,  or $L\ge 20|\Delta_0-\Delta_1|$, depending upon
the hypothesis in which the moments fall into as classified in
Theorem~\ref{secondtheorem}. The maximum required in
Problem~\ref{problem3} is 760 which is attained for $\alpha_0=1$,
$\alpha_1=-1$, $\alpha_2=1$, and $\alpha_3=-1$. We wonder if
alternating the signs of the moments and setting them
$\alpha_k=(-1)^k$ will always give the maximum in
Problem~\ref{problem3}.

For higher derivatives we can show that
$U:=120(3\alpha_1-12\alpha_2+10\alpha_3)$ and
$V:=120(\alpha_0-9\alpha_1+18\alpha_2-10\alpha_3)$ are in the
range of the second derivative. It does not seem to follow from
our Theorem~\ref{firsttheorem} applied to $f'$ that these values
are sharp, although the same idea of using a spline formed by a a
cubic and a quadratic may work.

 \vspace{0.2in}

\section{Higher values of $n$}

We have noted the following statement after the proof of
Proposition~\ref{enlargement}.

\begin{theorem}\label{somegeneraltheorem}
Given a continuously differentiable function satisfying the
Hausdorff moments constraints (\ref{thegiven}) ($n\ge 2$),  then
the range of the derivative contains the interval $[A_n,B_n]$,
where $A_n=\min \{D_k|k=1,2,...,n\}$ and $B_n=\max
\{D_k|k=1,2,...,n\}$, with $D_k$ given by (\ref{definitionOfD}).
Moreover, $[A_n,B_n]\subset [A_{n+1},B_{n+1}]$ for all $n\ge 2$.
\end{theorem}

\proof The first part follows with the same technique we have
employed over and over here using Theorem~\ref{mean}. For the
second part we are observing that the integrals which appear in
the denominators of (\ref{definitionOfD}),  are actually the well
known values of the beta function, i.e.
$B(k+1,n+2-k)=\int_0^1x^k(1-x)^{n+1-k}dx$. Using the established
formula for $B(\ ,\ )$,  we see that

$$B(k+1,n+2-k)=\frac{\Gamma(k+1)\Gamma(n+2-k)}{\Gamma(n+3)}=\frac{k!(n+1-k)!}{(n+2)!}=\frac{1}{(n+2){n+1 \choose k}}.$$

\n This gives us a new expression of the $D_{k,n}$ which is
basically in terms of the genuine Bernstein basis polynomials,
i.e. $b_{\nu,n}={n\choose \nu}x^{\nu}(1-x)^{n-\nu}$,
$\nu=0,1,...,n$:

\begin{equation}\label{definitionOfDkintermsofBernsteinpolynomials}
D_{k,n}:=-(n+2)\int_0^1f(x)\frac{d}{dx}(b_{k,n+1})dx,\ \
k=1,2,...,n.
\end{equation}

\n It is easy to check that
$x^k(1-x)^{n+1-k}+x^{k+1}(1-x)^{n-k}=x^k(1-x)^{n-k}$ which
basically gives the convex combination formula

$$D_{k,n}=\frac{n+2-k}{n+3}D_{k,n+1}+\frac{k+1}{n+3}D_{k+1,n+1},\ \ k=1,2,...,n,\ n\ge 1.$$

\vspace{0.2in}

\n These expressions imply the second claim of the theorem.\eproof
Let us observe that
(\ref{definitionOfDkintermsofBernsteinpolynomials}) implies the
following form for $D_{k,n}$
\begin{equation}\label{definitionOfDkintermsofBernsteinpolynomials}
D_{k,n}:=(n+1)(n+2)\int_0^1f(x)(b_{k,n}-b_{k-1,n})dx,\ \
k=1,2,...,n.
\end{equation}
\n which provides a simple way of computing  $D_k's$ in terms of
the moments $\alpha_0$, $\alpha_1$,..., $\alpha_n$. In what
follows we will describe yet another way of doing these
computations, and for that purpose we generalize first the
definitions of $\Delta_0$ and $\Delta_1$ in the following way

$$\Delta_k=\frac{(k+1)(k+2)^2(k+3)^2(k+4)}{2}\left|%
\begin{array}{ccc}
  \frac{1}{k+1} & \frac{1}{k+2} & \alpha_k \\
  \frac{1}{k+2} & \frac{1}{k+3} & \alpha_{k+1} \\
   \frac{1}{k+3} & \frac{1}{k+4} & \alpha_{k+2} \\
\end{array}%
\right|,k=0,1,2,....$$

\n or simply

$$\Delta_k=\frac{(k+3)(k+4)}{2}\alpha_{k+2}-(k+2)(k+3)\alpha_{k+1}+\frac{(k+1)(k+2)}{2}\alpha_k,\ \ k\ge 0. $$

\n There are some relations between the $D_k's$ and $\Delta_k's$
in general which we will include in the next proposition.

\begin{proposition}\label{deltadkrel} For $k\ge 0$ and $n\ge 2$, we have in general

\begin{equation}\label{connection}
\Delta_k=\frac{D_{k+2,k+2}-D_{k+1,k+1}}{2(k+3)}, k\ge 0.
\end{equation}

\begin{equation}\label{connection2}
D_{n,n}=6(2\alpha_1-\alpha_0)+ 2\sum_{k=0}^{n-2}(k+3)\Delta_k.
\end{equation}

\n Moreover, with the definitions of $A_n$ and $B_n$ from
Theorem~\ref{somegeneraltheorem}, $A_n=B_n$ if and only if
$\Delta_i=0$ for all $i=0,1,2,...,n-2$, if and only if there
exists a linear function with moments $\alpha_k$.
\end{proposition}

\proof  Let us observe that we can simply write

 $$x^k=x^{k+1}+x^k(1-x)=x^{k+1}+\frac{1}{k+1}b_{k,k+1}\Rightarrow$$
 $$k\alpha_{k-1}=(k+1)\alpha_{k}-\frac{1}{(k+1)(k+2)}D_{k,k}.$$

\n Using this last formula we can calculate the expression of
$\Delta_k$:

$$c\Delta_k=\left|%
\begin{array}{ccc}
  \frac{1}{k+1} & \frac{1}{k+2} & \frac{k+2}{k+1}\alpha_{k+1}-\frac{D_{k+1,k+1}}{(k+1)(k+2)(k+3)} \\
  \frac{1}{k+2} & \frac{1}{k+3} & \alpha_{k+1} \\
   \frac{1}{k+3} & \frac{1}{k+4} &  \frac{k+2}{k+3}\alpha_{k+1}+\frac{D_{k+2,k+2}}{(k+3)^2(k+4)}\\
\end{array}%
\right|=\left|%
\begin{array}{ccc}
  \frac{1}{k+1} & \frac{1}{k+2} & -\frac{D_{k+1,k+1}}{(k+1)(k+2)(k+3)} \\
  \frac{1}{k+2} & \frac{1}{k+3} &  0 \\
   \frac{1}{k+3} & \frac{1}{k+4} &  \frac{D_{k+2,k+2}}{(k+3)^2(k+4)}\\
\end{array}%
\right|,$$

\n where $c=\frac{2}{(k+1)(k+2)^2(k+3)^2(k+4)}$. This last
identity implies the formula (\ref{connection}). For the second
part of our statement we observe that $\Delta_i=0$ for all
$i=0,1,2,...n-2$ if and only if there exists a linear function $f$
with moments $\alpha_k$. In this case the range of the derivative
of $f$ consists of only one point and therefore by
Theorem~\ref{somegeneraltheorem} we must have $A_n=B_n$. For the
converse, again using Theorem~\ref{somegeneraltheorem} we obtain
that all $D_{i,j}$, $1\le i\le j$, $1\le j\le n$, have identical
values and so by (\ref{connection}) we get $\Delta_k=0$ for all
$k=0,1,...,n-2$.

Finally, let us observe that the equalities in (\ref{connection})
provide a telescopic sum for $D_{k,k}$ which allows one to arrive
at formula (\ref{connection2}).\eproof

The convexity relations can be used to calculate all the $D_k's$
from the $D_{k,k}$ and so formulae (\ref{connection2}) provide a
way of computing all the $D_k's$ in terms of determinants
$\Delta_i$.

For the case $\alpha_k=k+1$, $k=0,1,2,...,n$, we calculate $A_n$
and $B_n$ in a more precise way. This generalizes the problem in
\cite{dg}.

\begin{corollary}\label{generalexercise} Let $n\in \mathbb N$, $n\ge 2$, be fixed and $f$ be a continuously
differentiable satisfying (\ref{thegiven}) with $\alpha_k=k+1$,
$k=0,1,2,...,n$. Then, the values of $A_n$ and $B_n$ as defined in
Theorem~\ref{somegeneraltheorem} are
$$A_n=-n(n+1)(n+2),\ B_n=(n+1)(n+2)(2n+1).$$
\end{corollary}

\proof Using the formula for $\Delta_k$ we get

$$\Delta_k=\frac{(k+3)^2(k+4)}{2}-(k+2)^2(k+3)+\frac{(k+1)^2(k+2)}{2}=3k+7,\
\ k\ge 0. $$

Then using formula (\ref{connection2}) we obtain

$$D_{n,n}=18+2\sum_{k=0}^{n-2}(k+3)(3k+7)=(n+1)(n+2)(2n+1).$$

Now we can use the convexity relations and compute $D_{n-1,n}$:

$$D_{n-1,n}=\frac{1}{2}\left((n+2)D_{n-1,n-1}-nD_{n,n}\right)=-n(n+1)(n+2).$$

Next, if one calculates $D_{n-2,n}$, some surprise appears:

$$D_{n-2,n}=\frac{1}{3}\left((n+2)D_{n-2,n-1}-(n-1)D_{n-1,n}\right)=0.$$

\n Because of the convexity relation, it is easy to see that all
the other $D_{k,n}$, $k\le n-2$, are equal to zero. Therefore
$A_n=-n(n+1)(n+2)$ and $B_n=(n+1)(n+2)(2n+1)$. \eproof

Putting together what we did so far we now can say that for
$\alpha_k=k+1$, the bounds above are sharp if $n=2$ and the lower
bound is sharp if $n=3$.

\begin{theorem}\label{lowerboundofL}
For $n\ge 2$ fixed, with the definition of $A_n$ and $B_n$ as in
Theorem~\ref{somegeneraltheorem}, if $A_n<B_n$ it is not possible
to have $L=B_n-A_n$ in Problem~\ref{problem2}.
\end{theorem}

\proof By way of contradiction let us assume that $L=B_n-A_n$.
Hence, we can find a sequence of functions $f_m$, continuously
differentiable, such that
$Range(f'_m)\subset[A_n-\frac{1}{m},B_n+\frac{1}{m}]$ and
satisfying (\ref{thegiven}). Since $f_m'$ can be considered in
$L^2([0,1])$ we can find a subsequence of $f'_m$, say $f'_{m_k}$,
weakly convergent to a function $f\in L^2([0,1])$. This implies
that for every non-negative function $g\in L^2([0,1])$,
$$(A_n-\frac{1}{m})||g||_1 \le \int_0^1f_m'(x)g(x)dx\le
(B_n+\frac{1}{m})||g||_1,$$

\n where $||h||_1=\int_0^1h(x)dx$, $h\in L^1([0,1])$. Passing to
the limit as $m_k\to \infty$, we get

$$A_n||g||_1 \le \int_0^1f(x)g(x)dx\le
B_n||g||_1,\ g\in L^2([0,1]), g\ge 0.$$

\n This implies that $A_n\le f(x)\le B_n$ for a.e. $x\in [0,1]$,
by a standard measure theory argument. Since $B_n=D_k$ for some
$k=1,2,...,n$, and $f_m$ satisfy (\ref{thegiven}) we can say that

$$B_n=-\frac{\int_0^1f_m(x)\frac{d}{dx}[x^k(1-x)^{n+1-k}]dx}{|| x^k(1-x)^{n+1-k}||_1}=
\frac{\int_0^1f'_m(x) x^k(1-x)^{n+1-k}dx}{||
x^k(1-x)^{n+1-k}||_1}.$$

\n Letting $m_k\to \infty$ we obtain

$$B_n=
\frac{\int_0^1f(x) x^k(1-x)^{n+1-k}dx}{|| x^k(1-x)^{n+1-k}||_1}.$$

\n Hence, re-writing this yields
$$\int_0^1[B_n-f(x)]x^k(1-x)^{n+1-k}dx=0,$$ which in turn implies, by
what we have shown before about $f$, that $f(x)=B_n$ for a.e.
$x\in [0,1]$. Similarly, we arrive at the conclusion $f(x)=A_n$
for a.e. $x\in [0,1]$. Since we assumed $A_n<B_n$ we clearly get a
contradiction. Therefore, it remains that $L>B_n-A_n$.\eproof

This last theorem says that Problem~\ref{problem1} and
Problem~\ref{problem2} are completely different in nature. We must
admit that we do not have a definite answer to
Problem~\ref{problem2}, other than the trivial case $L=0$, in any
of the particular situations we have considered. We leave that to
the interested reader.


\begin{thebibliography}{9}
\bibitem{bartle} R. G. Bartle, {\it The elements of real
analysis}, Second Edition, 1976 John Wiley $\&$ Sons, Inc.
\bibitem{mch} Man-D. Choi, {\it Tricks or Treats with the Hilbert
Matrix},   The American Mathematical Monthly, Vol. 90, No. 5 (May,
1983), pp. 301-312.
\bibitem{gi2} G. Inglese {\it Christoffel functions and finite moment
problems}, Inverse Problems 11 (1995) 949-960.
\bibitem{gi3} G. Inglese {\it A note about the discretization of finite moment problems}, Inverse Problems 10 (1994) 401-414.
\bibitem{gi1} G. Inglese, {\it A note about minimum relative entropy solutions of finite moment
problems}, Numer. Funct. Anal.Optim. 16 (1995), no. 9-10,
1143–1153.
\bibitem{cezarlupu} C. Lupu, {\it Problem U37}, Mathematical Reflections no. 6 (2006).
\bibitem{rodriguez-seatzu} G. Rodriguez and S.  Seatzu {\it
On the solution of the finite moment problem.(English summary)} J.
Math. Anal. Appl. 171 (1992), no. 2, 321–333.
\bibitem{dvt} D. V. Thong, {\it The problem No 951}, The College Mathematics Journal, May (2011), 232-233. 
\end{thebibliography}
\end{document}